\documentclass[10pt,twoside]{article}
\usepackage{graphicx}
\usepackage{amsmath}
\usepackage{Latex-document}
\usepackage{amssymb}
\usepackage{amstext}
\usepackage{amscd}

\newtheorem{teo}{Th\'eor\`eme}[section]

\newtheorem{cor}[teo]{Corollaire}

\newtheorem{exe}[teo]{Exemple}

\newtheorem{ques}[teo]{Question}

\newcommand{\CC}{{\mathbb C}}
\newcommand{\RR}{{\mathbb R}}
\newcommand{\ZZ}{{\mathbb Z}}
\newcommand{\QQ}{{\mathbb Q}}
\newcommand{\NN}{{\mathbb N}}

\newcommand{\HH}{{\mathbb H}}

\newcommand{\GG}{{\mathbb G}}
\newcommand{\SSS}{{\mathbb S}}
\newcommand{\AAA}{{\mathbb A}}

\renewcommand{\thesection}{\arabic{section}}
\def\addsec{\addtocounter{section}{1} \setcounter{teo}{0}}

\title{\bf Th\'eorie Ergodique et G\'eom\'etrie \vskip -2mm Arithm\'etique\vskip 6mm}
\author{ Emmanuel Ullmo\vspace*{-0.5cm}\thanks{Universit\'e Paris-Sud
Orsay B\^at 425, 91405 Orsay Cedex France. E-mail:
ullmo@math.u-psud.fr}}
\date{\vskip -8mm}

\begin{document}

\maketitle

\thispagestyle{first} \setcounter{page}{197}

\begin{abstract}

\vskip 3mm

We will present several examples in which ideas from ergodic theory can be useful to study some problems in
arithmetic and algebraic geometry.

\vskip 4.5mm


\end{abstract}

\vskip 12mm

\section*{1. Introduction}\setzero \addsec
\vskip-5mm \hspace{5mm}

Le but de ce rapport est d'expliquer diff\'erentes techniques
permettant de montrer l'\'equidistribution de certains ensembles
de points de nature arithm\'etique sur des vari\'et\'es alg\'ebriques
d\'efinies sur des corps de nombres et de donner des applications
arithm\'etiques et g\'eom\'etriques de ces r\'esultats.

Si $X$ est une vari\'et\'e alg\'ebrique sur $\CC$ et $E$
une ensemble fini de $X(\CC)$ on note  $\vert E\vert$ son cardinal
et  $\Delta_{E}$ la mesure
de Dirac normalis\'ee
$$
\Delta_{E}=\frac{1}{\vert E\vert} \sum_{x\in E}\delta_{x}.
$$
Si $E_{n}$ est une suite d'ensembles finis
de $X(\CC)$ et $\mu$ une mesure de probabilit\'e sur $X(\CC)$,
on dit que les $E_{n}$ sont \'equidistribu\'es pour $\mu$ si pour
toute fonction continue born\'ee $f$ sur $X(\CC)$ on a
$$
\Delta_{E_{n}}(f)= \frac{1}{\vert E_{n}\vert}
\sum_{x\in E_{n}} f(x)\longrightarrow \int_{X(\CC)}f\mu.
$$

Soit $X$ une vari\'et\'e alg\'ebrique, une suite de points $x_n$
de $X$ est dite ``g\'en\'erique''  si pour toute sous-vari\'et\'e
$Y$ de $X$, $Y\neq X$, $\{n\in \NN\   ,x_n\in Y   \}$ est un
ensemble fini. (Il revient au m\^eme de dire que $x_n$ converge
vers le point g\'en\'erique pour la topologie de Zariski).

Andr\'e et Oort ont formul\'e un analogue de la conjecture de
Manin-Mumford d\'emontr\'ee par Raynaud \cite{Ra1} \cite{Ra2} dans le
cadre des vari\'et\'es de Shimura. Dans ces deux conjectures, on
dispose de points sp\'eciaux et de vari\'et\'es sp\'eciales.
Pour la conjecture de Manin-Mumford l'espace ambiant est une
vari\'et\'e ab\'elienne, les points sp\'eciaux sont les points de
torsion et les vari\'et\'es sp\'eciales sont les ``sous-vari\'et\'es
de torsion'' (translat\'es, par un point de torsion, d'une
sous-vari\'et\'e ab\'elienne). Pour la conjecture d'Andr\'e-Oort
l'espace ambiant est une vari\'et\'e de Shimura, les points sp\'eciaux
sont les points \`a multiplication complexe  (ou points CM)
et les sous-vari\'et\'es
sp\'eciales sont les ``sous-vari\'et\'es de type de Hodge''
(des composantes irr\'eductibles de translat\'es par un op\'erateur de
Hecke de sous-vari\'et\'es de Shimura). Nous pr\'eciserons ces
d\'efinitions plus bas. Dans les deux cas ces conjectures s'\'enoncent
sous la forme:  une composante irr\'eductible de l'adh\'erence de
Zariski d'un ensemble de points sp\'eciaux est une sous-vari\'et\'e
sp\'eciale.

Dans ce cadre
une suite de points $x_n$ de
$X$ ($X$ vari\'et\'e ab\'elienne ou $X$ vari\'et\'e de Shimura)
est dite ``stricte''  si pour toute sous-vari\'et\'e sp\'eciale  $Y$
de $X$,  $Y\neq X$,    $\{n\in \NN\   ,x_n\in Y   \}$ est un ensemble fini.
On remarque qu'avec ces d\'efinitions les conjectures d'Andr\'e-Oort
et de Manin-Mumford se retraduisent de la mani\`ere suivante:
Toute suite stricte de points sp\'eciaux est g\'en\'erique.

Une cons\'equence g\'eom\'etrique
(conjecturale pour les vari\'et\'es de Shimura)
que l'on obtient en consid\'erant
l'adh\'erence de Zariski de l'ensemble des points sp\'eciaux d'une
sous-vari\'et\'e $M$ de $X$ est l'existence d'un ensemble fini
$\{S_{1},\ldots,S_{r}\}$ de sous-vari\'et\'es sp\'eciales avec
$S_{i}\subset M$ telle que toute sous-vari\'et\'e sp\'eciale $S\subset
M$ est
contenue
dans l'un des $S_{i}$.

Dans la premi\`ere partie nous d\'ecrivons des r\'esultats
d'\'equidistribution pour des suites de points de petite hauteur
sur des vari\'et\'es alg\'ebriques
utilisant la g\'eom\'etrie d'Arakelov. Le r\'esultat le plus
marquant est la r\'esolution de la conjecture de Bogomolov
(qui g\'en\'eralise la conjecture de Manin-Mumford
et en donne une nouvelle d\'emonstration)
pour les vari\'et\'es ab\'eliennes due \`a Zhang \cite{Zh}
et \`a l'auteur
du rapport \cite{Ul}.

Dans la deuxi\`eme partie
nous expliquons des r\'esutats d'\'equidistribution
de points  de Hecke sur des vari\'et\'es de la forme
$X=\Gamma\backslash G(\RR)$ pour un groupe alg\'ebrique
semi-simple et simplement connexe $G$ et un r\'eseau $\Gamma$.
Les m\'ethodes combinent th\'eorie spectrale et th\'eorie des
repr\'esentations.

Dans la troisi\`eme partie nous pr\'esentons
des \'enonc\'es largement conjecturaux pour l'\'equidistribution des
points \`a multiplication complexe des vari\'et\'es
de Shimura. La th\'eorie analytique des nombres
via les familles de fonctions $L$ et la th\'eorie
des formes automorphes y jouent un r\^ole central.

Dans une derni\`ere
partie nous expliquons comment la th\'eorie de Ratner et Margulis
permet de d\'emontrer des r\'esultats d'\'equidistribution pour
des suites de sous-vari\'et\'es ``fortement sp\'eciales'' (appartenant
\`a une classe assez large de sous-vari\'et\'es sp\'eciales
de dimension positive)
des vari\'et\'es de Shimura. Nous expliquerons la relation avec
 la cons\'equence g\'eom\'etrique de
la conjecture d'Andr\'e-Oort pr\'ec\'edemment d\'ecrite.

\section*{2. Equidistribution des points de petite hauteur} \addsec \setzero

\vskip-5mm \hspace{5mm}

\begin{exe}\label{exe1}\rm
On prend $X=\GG_{m}$, $E_{n}$ l'ensemble des racines $n$-i\`eme de
l'unit\'e, $E_{n}$ est \'equidistribu\'e pour la mesure uniforme
sur le cercle unit\'e $\frac{d\alpha}{2\pi}$. En utilisant
l'irr\'eductibilit\'e du polynome cyclotomique on voit que l'orbite
sous Galois d'une racine $n$-i\`eme primitive de l'unit\'e est aussi
\'equidistribu\'ee pour $\frac{d\alpha}{2\pi}$.
\end{exe}

\begin{exe}\rm
 On prend $X=E$ une courbe elliptique sur $\CC$ et $E_{n}$ l'ensemble
 des points de $n$ torsion, alors $E_{n}$ est \'equidistribu\'e pour
 la mesure de Harr normalis\'ee sur $E(\CC)$. Si $E$ est d\'efini
 sur un coprs de nombres $K$ et $E$ n'a pas de multiplication complexe,
 par le th\'eor\`eme de l'image ouverte de Serre, pour tout nombre
 premier $p$ assez grand le groupe de Galois agit transitivement sur
 les points d'ordre $p$. On en d\'eduit encore que les orbites sous
 Galois des points d'ordre $p$ sont \'equidistribu\'ees pour la mesure
 de Haar normalis\'ee.
 \end{exe}

La th\'eorie d'Arakelov a permis de comprendre  ces \'enonc\'es
d'une mani\`ere bien plus g\'en\'erale. On montre \cite{SUZ}
pour une vari\'et\'e
arithm\'etique un th\'eor\`eme g\'en\'eral d'\'equidistribution
des orbites sous Galois de suite g\'en\'eriques de points dont la hauteur (\`a
la Arakelov) tend vers $0$. Les exemples pr\'ec\'edents correspondent
\`a des suites de points de hauteurs nulles. Pour les vari\'et\'es
ab\'eliennes on obtient avec Szpiro et Zhang
le r\'esultat suivant (qui donne des
informations
nouvelles m\^eme pour les points de torsion des courbes elliptiques
\`a multiplication complexe):
\begin{teo}\label{teo1}{\rm \cite{SUZ}}
Soit $A$ une vari\'et\'e ab\'elienne sur un corps de nombres $K$. On
note $h_{NT}$ la hauteur de N\'eron-Tate sur les points alg\'ebriques
de $A$ (associ\'ee \`a un fibr\'e inversible ample sym\'etrique sur
$X$). Soit $x_n$ une suite
 g\'en\'erique  de points alg\'ebriques de $A$ telle que
 $h_{NT}(x_n)$ tend vers $0$. Pour toute place \`a l'infini $\sigma$
 l'orbite sous Galois de $x_{n}$ est \'equidistribu\'ee pour la mesure
 de Haar normalis\'ee $d\mu_{\sigma}$ de $A_{\sigma}(\CC)$.
\end{teo}

 L'analogue de cet \'enonc\'e pour $\GG_{m}^r$ a \'et\'e montr\'e par
 Bilu \cite{Bi} sans th\'eorie d'Arakelov. Une extension
 pour certaines vari\'et\'es semi-ab\'eliennes de ces r\'esultats a
 \'et\'e obtenue par Chambert-Loir \cite{Ch}
 par des m\'ethodes Arakeloviennes.
 On peut aussi comprendre gr\^ace aux travaux de Autissier \cite{Au}
 l'exemple \ref{exe1} comme un cas particulier
 de th\'eor\`eme d'\'equidistribution vers la mesure d'\'equilibre
 d'un compact de capacit\'e $1$ de l'orbite sous Galois
 d'une suite de points entiers alg\'ebriques.

 On trouvera dans \cite{Zh2} comment on obtient la conjecture
 de Bogomolov en produisant une contradiction sur les mesures
 limites de suites de mesures associ\'ees \`a des orbites sous
 Galois de points de petite hauteur.
 Retenons l'\'enonc\'e suivant
 d\^u \`a l'auteur \cite{Ul} pour les courbes de genre $g\ge 2$ dans leur
 jacobienne et \'etendu en dimension arbitraire par Zhang \cite{Zh}:

 \begin{teo}\label{teo2}
Soit $X$ une sous-vari\'et\'e
 d'une vari\'et\'e ab\'elienne $A$ d\'efinie sur un corps de nombres $K$.
 Gr\^ace \`a la conjecture de Manin-Mumford d\'emontr\'ee par Raynaud
 \cite{Ra2}, on sait qu'il existe des sous-vari\'et\'es de torsion
 (\'eventuellement de dimension $0$)
 $\{T_{1},\ldots,T_{r}\}$, $T_{i}\subset X$ tels que si $T\subset X$ est une
 sous-vari\'et\'e de torsion alors $T\subset T_{i}$ pour un certain
 $i$. Il existe alors $c>0$ tel que si $P$ est un point alg\'ebrique
 de $X$ et $P\notin \cup_{i=1}^r T_{i}$ alors $h_{NT}(P)\ge c$.
\end{teo}

\section*{3. Equidistribution des points de Hecke} \addsec \setzero

\vskip-5mm \hspace{5mm}

Soient $G$ un groupe alg\'ebrique lin\'eaire presque simple et simplement
connexe sur $\QQ$, $\Gamma\subset G(\QQ)$ un r\'eseau de congruence
et $X=\Gamma\backslash G(\RR)$. Soit $\mu_{0}$ la mesure invariante
normalis\'ee sur $X$.
Pour tout $a\in G(\QQ)$ on a
une d\'ecomposition
$$
\Gamma a\Gamma=\cup_{i=1}^{deg(a)}\Gamma a_{i}
$$
avec $deg(a)=\vert \Gamma \backslash \Gamma a\Gamma\vert\in \NN$. Pour
tout $x\in X$,
on note $T_{a}.x$ l'ensemble des $a_{i}x$ compt\'e avec multiplicit\'e.
L'op\'erateur de Hecke $T_{a}$ ainsi d\'efini est une correspondance
de degr\'e $deg(a)$ sur $X$; il induit une op\'eration sur
les espaces de fonctions $L^2(X,\mu_{0})$
(fonctions de carr\'es int\'egrables sur $X$) et $C^0_{b}(X)$
(fonctions continues born\'ees sur $X$) par
$$
T_{a}.f(x)=\frac{1}{deg(a)}\sum_{i=1}^{deg(a)} f(a_{i}x).
$$
Avec Clozel et Oh nous obtenons \cite{COU}:
\begin{teo} On suppose que le $\QQ$-rang de $G$ est diff\'erent de $0$.
Soit $a_{n}\in G(\QQ)$ une suite telle que $deg(a_{n})\rightarrow
\infty$. Pour tout $x\in X$ les $T_{a_{n}}.x$ sont
\'equidistribu\'es pour $\mu_{0}$.
De plus pour tout $f\in L^2(X,\mu_{0})$ on a la convergence $L^2$
$$
\Vert T_{a_{n}}f-\int_{X} f.\mu_{0}\Vert_{L^2}\longrightarrow 0.
$$
\end{teo}

On a en fait des r\'esultats aussi dans le cas ou le $\QQ$-rang de
$G$ vaut $0$. La m\'ethode de d\'emonstration fournit des
estimations tr\`es pr\'ecises pour la vitesse de convergence dans
le th\'eor\`eme $L^2$. Si on dispose de plus de r\'egularit\'e sur
$f$ (par exemple $f$ $C^{\infty}$ \`a support compact), cete
vitesse est obtenue aussi pour la convergence simple (ou uniforme
sur les compacts). Pour $G=SL_{n}$ ($n\ge 3$) ou $G=Sp_{2n}$
($n\ge 2$) ces estimations sont essentiellement optimales.

On montre par des m\'ethodes classiques que l'\'enonc\'e
de convergence simple du th\'eor\`eme se d\'eduit de l'\'enonc\'e
$L^2$.
Pour montrer le th\'eor\`eme $L^2$ on \'ecrit la d\'ecomposition
spectrale de $L^2(X,\mu_{0})$ sous la forme ad\'elique.
Une  fonction $\phi$ intervenant dans la d\'ecomposition spectrale
est alors propre pour les op\'erateurs de Hecke et
les valeurs propres s'interpr\`etent comme des coefficients
matriciaux de  repr\'esentations locales associ\'ees \`a $\phi$.
Pour montrer le th\'eor\`eme sous la forme $L^2$, on
doit montrer que $T_{a_{n}} \phi \rightarrow 0$ quand $n\rightarrow
\infty$ au sens $L^2$.
On se  ram\`ene ainsi \`a contr\^oler
la d\'ecroissance de ces coefficients matriciaux. En $\QQ$-rang $r\ge 2$
on dispose d'assez d'informations sur le dual unitaire pour
conclure gr\^ace aux travaux de Oh (\cite{Oh}, th\'eor\`eme 5.7).
En $\QQ$-rang $1$
on utilise un principe de restriction \`a la Burger-Sarnak
en une place finie d\'emontr\'e dans \cite{CU} et une approximation
de la conjecture de Ramanujan pour $SL_{2}$.

\section*{4. Equidistribution des points CM des vari\'et\'es de Shimura} \addsec \setzero

\vskip-5mm \hspace{5mm}

Nous devons
pr\'eciser un peu les
d\'efinitions relatives aux vari\'et\'es de Shimura afin d'expliquer
ce que l'on entend par l'\'equidistribution des points CM.

Soit  $(G,X)$ une donn\'ee de Shimura; $G$ est un groupe
alg\'ebrique r\'eductif sur $\QQ$ et $X$ est une $G(\RR)$
classe de conjuguaison de morphismes
$$
h:\ \SSS \longrightarrow G_{\RR}
$$
($\SSS=\mbox{Res }_{\CC/\RR} \mathbb{G}_m$ est le tore de Deligne)
v\'erifiant les 3 propri\'et\'es de Deligne \cite{De1}
\cite{De2}. Les composantes irr\'eductibles de $X$ sont
alors des domaines sym\'etriques hermitiens.

Soient $\AAA_f$ l'anneau des ad\`eles finies de $\QQ$ et
$K$ un sous-groupe compact ouvert de $G(\AAA)$, on d\'efinit sur
le corps $\CC$ la vari\'et\'e de Shimura
$$
Sh_K(G,X)=G(\QQ)\backslash {X\times G(\AAA_f)/K}.
$$
On v\'erifie que $Sh_K(G,X)$ est une r\'eunion finie de quotients de composantes \linebreak irr\'eductibles de $X$
par des sous-groupes de congruences de $G(\QQ)$. Par ailleurs $Sh_K(G,X)$ a un ``mod\`ele canonique'' sur un corps
de nombres $E(G,X)$ ne d\'ependant que de la donn\'ee de Shimura $(G,X)$.

Soit $(G_1, X_1)$ une sous-donn\'ee de Shimura de $(G,X)$, on dispose
alors d'une application canonique
$$
f:\ Sh_{K\cap G_1(\AAA_f)}\longrightarrow Sh_K(G,X).
$$
Une sous-vari\'et\'e de type de Hodge est une composante irr\'eductible
d'un translat\'e de l'image d'un tel morphisme par une correspondance
de Hecke.  (Moonen \cite{Mo} caract\'erise ces sous-vari\'et\'es en
termes
de variations de structures de Hodge, d'o\`u le nom.)

Pour $h: \SSS \rightarrow G_{\RR}$, $h\in X$, on d\'efinit le groupe
de Mumford-Tate $MT(h)$ de $h$ comme le plus petit
$\QQ$-sous-groupe $H$ de $G$ tel que $h$ se factorise
par $H_{\RR}$. Si $MT(h)$ est un tore, on dit que $h$ est
sp\'ecial. Les points sp\'eciaux de $Sh_K(G,X)$ sont les points de la
forme $[h,gK]$ avec $g\in G(\AAA_f)$ et $h$ sp\'ecial.

Fixons $h_0\in X$ un \'el\'ement sp\'ecial et $T_0=MT(h_0)$.
L'ensemble
$$
S(h_0)=\{[h_0,gK],\ \  g\in G(\AAA_f)   \}
$$
est appel\'e ensemble des points sp\'eciaux de ``type $h_0$'' de $X$.
On a une action de $T_0(\AAA_f)$ sur $S(h_0)$ donn\'ee par
$t.[h_0,gK]=[h_0,tgK]$. Pour tout $g\in G(\AAA_f)$, l'orbite sous
$T_0(\AAA_f)$ de $[h_0,gK]$ est finie, on appelle ``orbite torique''
de $[h_0,gK]$ cette orbite. La premi\`ere question naturelle
est

\begin{ques}\label{ques1}\rm
Soit $x_{n}=[h_{n},g_{n}K]$ une suite g\'en\'erique de points sp\'eciaux
de  $S=Sh_K(G,X)$. Est-il vrai que l'orbite torique de $x_{n}$
est \'equidistribu\'ee pour la mesure invariante normali\'see de
$Sh_K(G,X)$.
\end{ques}

Notons qu'il n'est d\'ej\`a pas \`a priori \'evident de pr\'evoir
la proportion des points de l'orbite torique dans les composantes
de $S$. Il est peut-\^etre plus r\'ealiste de travailler dans
chaque composante connexe de $S$ (comme dans la derni\`ere partie
de ce texte). Nous tairons dans la suite ces probl\`emes de non
connexit\'e.

Les premiers r\'esultats pour ces questions sont dus \`a Duke \cite{Du}
pour la courbe modulaire $Y(1)=SL(2,\ZZ)\backslash \HH$. Il montre
l'\'equidistribution des points \`a multiplication complexe par
l'anneau des entiers $O_K$ quand le discriminant tend vers l'infini.
Nous expliquons dans \cite{CU}, en utilisant
en plus des r\'esultats sur l'\'equidistribution
des points de Hecke comment obtenir
l'\'equidistribution des points \`a multiplication complexe
par un ordre arbitraire de $O_K$ quand le discriminant tend vers
l'infini. Nous pensons plus g\'en\'eralement que la  question
\ref{ques1} est li\'ee aux probl\`emes d'\'equidistribution des points
de Hecke d\'ecrits pr\'ec\'edemment.

Des r\'esultats pour l'\'equidistribution des orbites toriques
de points CM sont annonc\'es par S. Zhang \cite{Zh3}
pour les courbes de Shimura
et plus g\'en\'eralement des vari\'et\'es de Shimura
de type quaternionique via un avatar de la formule de Gross-Zagier.
Pour les vari\'et\'es modulaires de Hilbert
des r\'esultats de ce type sont annonc\'es ind\'ependamment
par P. Cohen \cite{Coh} (par la m\'ethode originale de Duke) et
par S. Zhang.

Les m\'ethodes pour prouver ces \'enonc\'es comportent trois
\'etapes que l'on va d\'ecrire de mani\`ere impr\'ecise pour
la concision de ce rapport. Soit $S$ une vari\'et\'e de Shimura,
soit $f$ une fonction non constante
intervenant dans la d\'ecomposition spectrale
de $S$, soit $x_n\in S$ une suite de points CM et $E_n$
son orbite torique. On doit montrer que
\begin{equation}\label{1}
\lim_{n\rightarrow \infty}\frac{1}{\vert E_n\vert}\sum_{y\in E_n}
f(y)=\int_S f d\mu_0.
\end{equation}
La fonction $f$ est alors une forme automorphe. La premi\`ere \'etape
est de montrer une ``formule de classe'' reliant
$\frac{1}{\vert E_n\vert}\sum_{y\in E_n}
f(y)$
\`a la valeur  de  la fonction $L$ de $f$, tordue
par une forme automorphe que l'on d\'efinit \`a partir de $E_n$,
au point critique. Ce type de formule est obtenu par Waldspurger
\cite{Wa} pour des alg\`ebres de quaternions sur un corps de nombres
$F$ et revisit\'e par Zhang \cite{Zh3}dans le but d'obtenir les r\'esultats
d'\'equidistribution.

Une fois la formule de classe \'etablie, on dispose d'une famille
de fonctions $L$ index\'ee par les entiers. On d\'efinit \`a
partir de l'\'equation fonctionnelle de ces fonctions une notion
de ``conducteur analytique'' $q_n$. L'hypoth\`ese de Riemann (ou
de Lindel\"of) pr\'evoit une borne en $0(q_n^{\epsilon})$ pour la
valeur critique de la fonction $L$ consid\'er\'ee. Dans tous les
exemples consid\'er\'es, il est remarquable que pour montrer
l'\'equidistribution il faut am\'eliorer la borne triviale
(donn\'ee par le principe de  convexit\'e de Phragmen-Lindel\"of).
Ce genre de questions a re\c cu une attention consid\'erable en
th\'eorie analytique des nombres et a \'et\'e r\'esolue dans de
nombreux cas. On pourra consulter  la s\'erie de papiers
\cite{DFI} et \cite{Fr} pour une pr\'esentation des principaux
r\'esultats et applications de ce cercle d'id\'ees. Notons que la
d\'emonstration de l'\'equidistribution des orbites toriques de
points CM sur  les vari\'et\'es modulaires de Hilbert utilise les
r\'esultats spectaculaires r\'ecents \cite{CPS}.

Pour les applications \'eventuelles \`a des \'enonc\'es
arithm\'etiques, il para\^\i t important de remplacer les orbites
toriques par les orbites sous Galois. De mani\`ere g\'en\'erale si
$[h,gK]$ est un point CM d'une
vari\'et\'e de Shimura, $T=MT(h)$  est le tore associ\'e
et $E=E(T,h)$ est le corps reflexe de la vari\'et\'e de Shimura
associ\'e \`a la donn\'ee de Shimura $(T,h)$,
l'action de Galois  (cf \cite{De1}, \cite{De2} )se factorise \`a travers l'action de $T(\AAA_f)$ via un
morphisme de r\'eciprocit\'e (et la th\'eorie du corps de classe).
$$
r :\mbox{Res}_{E/\QQ} \mathbb{G}_{m,E}\longrightarrow T
$$
qui induit un morphisme non surjectif en g\'en\'eral
$$
r :\mbox{Res}_{E/\QQ} \mathbb{G}_{m,E}(\AAA_f)\longrightarrow
T(\AAA_f).
$$
On s'attend n\'eanmoins \`a une r\'eponse positive \`a la question
suivante:
\begin{ques} \rm
Soit $x_n$ une suite g\'en\'erique de points CM sur une vari\'et\'e de Shimura $S$, est-il vrai que les orbites
sous Galois $O(x_n)$ sont \'equidistribu\'ees dans $S$ pour la mesure invariante?
\end{ques}

De mani\`ere encore plus optimiste, on esp\`ere
(par analogie avec le cas des vari\'et\'es ab\'eliennes)
que le m\^eme r\'esultat est encore vrai pour des suites
strictes de points CM. Ce serait une cons\'equence de
la conjecture d'Andr\'e-Oort et de la question pr\'ec\'edente.
Notons que nous esp\'erons que des r\'esultats d'\'equidistribution
pour les points CM soient en fait une \'etape pour montrer la
conjecture en question . (C'est au moins ce qui se passe dans
le cas des vari\'et\'es ab\'eliennes).

\section*{5. Equidistribution de sous-vari\'et\'es sp\'eciales} \addsec \setzero

\vskip-5mm \hspace{5mm}

Cette partie d\'ecrit un travail \cite{CU2} en cours de pr\'eparation
en commun avec L. Clozel.
Soit $S$ une composante irr\'eductible d'une vari\'et\'e de Shimura.
Une cons\'equence
g\'eom\'etrique frappante de la conjecture d'Andr\'e et Oort
est la suivante:
Soit $Y$ une sous-vari\'et\'e de $S$, il existe un ensemble
fini $\{S_{1},\ldots,S_{r}\}$ de sous-vari\'et\'es sp\'eciales
avec $S_{i}\subset Y$ pour tout $i$ tel que toute vari\'et\'e
sp\'eciale $Z$ de $S$ contenue dans $Y$ est en fait contenue dans
un des $S_{i}$.

Supposons que $S$ est une composante irr\'eductible de $Sh_K(G,X)$
pour un groupe $G$ que l'on suppose adjoint (pour simplifier).
On a vu qu'une   sous-vari\'et\'e sp\'eciale $M$ est associ\'ee
\`a une sous-donn\'ee de Shimura $(G_{1},X_{1})$. Si $G_{1}$
est semi-simple et $X_{1}$ contient un point sp\'ecial $x_{1}$
tel que le tore associ\'e $T=MT(x_{1})\subset G_{1}$ est tel que
$T_{\RR}$ est
un tore maximal compact de $G$, on dit que $M$ est fortement sp\'eciale.
Par exemple les vari\'et\'es modulaires de Hilbert
(associ\'ees \`a des corps totalement r\'eels de degr\'e $n$
sur $\QQ$)
sont fortement
sp\'eciales dans l'espace de module  ${\cal A}_{n}$
 des vari\'et\'es ab\'eliennes
principalement polaris\'ees de dimension $n$. On peut montrer:

\begin{teo} \label{teo1.1}
     Soit $Y$ une sous-vari\'et\'e d'une vari\'et\'e de
Shimura $S$.
Il existe un
   ensemble fini $\{S_{1},\ldots,S_{k}\}$
de sous-vari\'et\'es fortement sp\'eciales de dimension positive
$S_{i}\subset
Y$ tel que si $Z$ est une sous-vari\'et\'e fortement sp\'eciale de dimension
positive avec $Z\subset Y$ alors $Z\subset S_{i}$ pour un certain
$i\in \{1,\dots,k\}$.
\end{teo}

Notons que cet \'enonc\'e ne dit rien sur les sous-vari\'et\'es
sp\'eciales de dimension $0$ (les points sp\'eciaux), notons
cependant le corollaire suivant:

\begin{cor}\label{cor1}
 Soit $Y$ une sous-vari\'et\'e stricte de ${\cal A}_{n}   $,
 il existe au plus un nombre fini de sous-vari\'et\'es modulaires de
 Hilbert contenu dans $Y$.
\end{cor}

Le th\'eor\`eme \ref{teo1.1} se d\'eduit d'un \'enonc\'e ergodique.
Toute sous-vari\'et\'e sp\'eciale $Z$ de $S$ est muni
d'une mani\`ere canonique d'une mesure de probabilit\'e $\mu_{Z}$.

\begin{teo}\label{teo2.2}
Soit $S_{n}$ une suite de sous-vari\'et\'es fortement sp\'eciales
Soit $\mu_{n}$ la mesure de probabilit\'e associ\'ee \`a $S_{n}$.
Il existe une sous-vari\'et\'e fortement sp\'eciale $Z$ et une
sous-suite $\mu_{n_{k}}$ qui converge faiblement vers $\mu_{Z}$.
De plus $Z$ contient $S_{n_{k}}$ pour tout $k$ assez grand.
\end{teo}

On obtient la preuve du th\'eor\`eme \ref{teo1.1} en consid\'erant
une suite de sous-vari\'et\'es fortement sp\'eciales maximales
$S_{n}$ parmi les sous-vari\'et\'es fortement sp\'eciales
contenues dans $Y$. En passant \`a une sous-suite on peut supposer
que $\mu_{n}$ converge faiblement vers $\mu_{Z}$. Comme le support
de $\mu_{Z}$ est contenu dans $Y$, on en d\'eduit que $Z\subset
Y$. Par la maximalit\'e des $S_{n}$ et le fait que $S_{n}\subset
Z$ pour tout $n$ assez grand, on en d\'eduit que la suite $S_{n}$
est stationaire.

On peut aussi r\'e\'ecrire cet \'enonc\'e avec la terminologie
de \cite{SUZ}. On dit qu'une suite $S_{n}$ de sous-vari\'et\'es fortement
sp\'eciales est stricte si pour toute sous-vari\'et\'e fortement
sp\'eciale $M$ de $S$,
$$
\{n\in \NN, \mbox{ } S_{n}\subset M\}
$$
est fini. On peut d'ailleurs prendre dans cette d\'efinition
$M$ sp\'eciale car une sous-vari\'et\'e sp\'eciale contenant une
sous-vari\'et\'e fortement sp\'eciale est automatiquement fortement
sp\'eciale. Dans ce language le th\'eor\`eme \ref{teo2.2}
admet comme corollaire imm\'ediat:

\begin{cor}
Soit $S_{n}$ une suite stricte de sous-vari\'et\'es
fortement sp\'eciales de $S$. Soit $\mu_{n}$ et $\mu$ les
mesures de probabilit\'es associ\'ees sur $S_{n}$ et $S$.
La suite $\mu_{n}$ converge faiblement vers $\mu$.
\end{cor}

On peut appliquer cet \'enonc\'e \`a des suites de
sous-vari\'et\'es fortement sp\'eciales maximales. La condition
d'\^etre stricte signifie alors de ne pas avoir de sous-suites
constantes. C'est par exemple le cas pour les vari\'et\'es
modulaires de Hilbert dans le modules des vari\'et\'es
ab\'eliennes principalement polaris\'ees ${\cal A}_{n}$.

La preuve du th\'eor\`eme \ref{teo2.2} repose sur des r\'esulats de
Mozes et Shah \cite{MoSh} qui pr\'ecisent la conjecture de Raghunathan
d\'emontr\'ee par Ratner \cite{Rat1}. Si $S=\Gamma\backslash
G(\RR)/K_{\infty}$ pour un sous-groupe compact maximal $K_{\infty}$
et un r\'eseau de congruence $\Gamma$, on note
$\Gamma^+=G(\RR)^+\cap \Gamma$ et $\tilde{S}=\Gamma^+\backslash
G(\RR)^+$. Si $H$ est un sous-groupe semi-simple de $G(\RR)^+$
tel que $\Gamma^+\cap H$
est un r\'eseau de $H$ alors $M_{H}=\Gamma^+\cap H\backslash H$ est
ferm\'e dans  $\tilde{S}$ et est muni canoniquement d'une mesure de
probabilit\'e
$H$-invariante $\mu_{H}$.

Si $M_{H_{n}}$ est une suite de telles sous-vari\'et\'es de
$\tilde{S}$, le th\'eor\`eme de Mozes Shah \cite{MoSh} permet sous
certaines hypoth\`eses; au besoin en passant \`a une sous-suite;
de montrer la convergence faible de $\mu_{H_{n}}$ vers une
mesure $\mu_{H}$ canoniquement associ\'ee
\`a un $M_{H}=\Gamma^+\cap H\backslash H$. En g\'en\'eral
les sous-groupes $H_{n}$ n'induisent pas de sous-vari\'et\'es
sp\'eciales sur $S$ car $H_{n}$ n'est pas toujours r\'eductif et
m\^eme si $H_{n}$ est r\'eductif
l'espace sym\'etrique associ\'e \`a
$H_{n}$ n'a aucune raison d'\^etre hermitien. Un des points
clefs de la d\'emonstration est de v\'erifier que si les $H_{n}$
induisent des sous-vari\'et\'es fortement sp\'eciales il en est de
m\^eme pour $H$. Pour passer de r\'esultats sur $\tilde{S}$
\`a des r\'esultats sur $S$ on utilise aussi des r\'esultats
de Dani et Margulis (\cite{DaMa}thm. 2) qui donnent des crit\`eres
de retour vers des compacts pour des flots unipotents sur $\tilde{S}$.

\label{lastpage}

\end{document}